Konrad-Zuse-Zentrum für Informationstechnik Berlin
Takustr. 7, D-14195 Berlin - Dahlem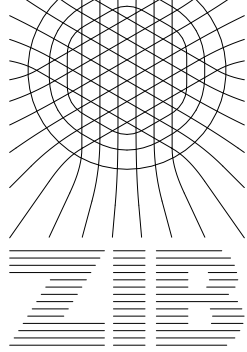

Wolfram Koepf

# On a Problem of Koornwinder

Preprint SC 96–52 (December 1996)

# On a Problem of Koornwinder

Wolfram Koepf

koepf@zib.de

**Abstract:**

In this note we solve a problem about the rational representability of hypergeometric terms which represent hypergeometric sums. This problem was proposed by Koornwinder in [4].

## 1 Hypergeometric Functions and Zeilberger's Algorithm

Zeilberger's algorithm ([6]–[7], see also [3], [1]) determines recurrence equations for hypergeometric functions

$$S(n) := {}_pF_q\left(\begin{array}{cccc}\alpha_1 & \alpha_2 & \cdots & \alpha_p \\ \beta_1 & \beta_2 & \cdots & \beta_q\end{array}\bigg| x\right) = \sum_{k=0}^{\infty} A_k\, x^k = \sum_{k=0}^{\infty} \frac{(\alpha_1)_k \cdot (\alpha_2)_k \cdots (\alpha_p)_k}{(\beta_1)_k \cdot (\beta_2)_k \cdots (\beta_q)_k} \frac{x^k}{k!}$$

whose upper parameters $\alpha_k$ and lower parameters $\beta_k$ are rational-linear in a variable $n$, whenever the term ratio

$$\frac{A_{k+1}}{A_k} = \frac{(k+\alpha_1)\cdot(k+\alpha_2)\cdots(k+\alpha_p)}{(k+\beta_1)\cdot(k+\beta_2)\cdots(k+\beta_q)\cdot(k+1)} \in \mathbb{Q}(k,n)$$

is a rational function in both $n$ and $k$. As usual $(a)_k = a(a+1)\cdots(a+k-1)$ denotes the *Pochhammer symbol*. We call the summand $A_k x^k$ a *hypergeometric term*. The resulting recurrence equation has polynomial coefficients with respect to $n$. If it is of first order, the sum has a rational term ratio with respect to $n$, and hence is itself a hypergeometric term.

In [4] Koornwider asked the question whether an application of Zeilberger's algorithm might generate a hypergeometric term whose upper and lower parameters are not rational assuming the parameters of the input summand are rational:

> **Problem 6.1.** If Zeilberger's algorithm succeeds, can $S(n)/S(n-1)$ then always be factorized as a quotient of products of linear forms over $\mathbb{Z}$ in $n$ and the parameters?

In this note, we will answer Koornwinder's question in the negative, by providing several counterexamples.

Note that Koornwinder's question in principle is independent of Zeilberger's algorithm, and asks whether there are hypergeometric sums that can be represented by hypergeometric terms with nonrational parameters. None example of this type can be found in the literature, see in particular the rather extensive mathematical dictionary on hypergeometric function identities [5]. Nevertheless we will use Zeilberger's algorithm to find our counterexamples.

In [1], we presented a Maple implementation for an extension of Zeilberger's algorithm, which is available through the `sumtools` package of Maple V.4 (or through the share library package `summation`). In this note we use a new implementation, developed in the book [2], which can be obtained from `http://www.zib.de/koepf/code` and `http://www.zib.de/koepf/simpcomb`,

but the recurrence equation calculations can also be done with the `sumtools` package of Maple V.4.

In [1], we gave the following example:

$$_2F_1\left(\begin{array}{c}1/2, -2n \\ 2n+3/2\end{array} \bigg| 3 \pm 2\sqrt{2}\right) = \sum_{k=0}^{2n}(-1)^k \frac{\binom{2n}{k}\binom{2n+k+1}{k}}{\binom{4n+2k+2}{2k}}\left(3 \pm 2\sqrt{2}\right)^k = \frac{(3/4)_n \, (5/4)_n}{(7/8)_n \, (9/8)_n}.$$

This example shows that despite the radicals of the input series, the output is rational. More examples of this type can be found in [5]. The above example can be generated by the Maple command (with the package `code`, or in Maple V.4, after `with(sumtools)`: `readlib('sum/simpcomb'):`)

```
> sumrecursion(hyperterm([1/2,-2*n],[2*n+3/2],3+2*sqrt(2),k),k,S(n));
```
$$-(\,8\,n+9\,)\,(\,8\,n+7\,)\,\mathrm{S}(\,n+1\,) + 4\,\mathrm{S}(\,n\,)\,(\,4\,n+5\,)\,(\,4\,n+3\,) = 0$$

or by

```
> sumrecursion((-1)^k*binomial(2*n,k)*binomial(2*n+k+1,k)/
> binomial(4*n+2*k+2,2*k)*(3+2*sqrt(2))^k,k,S(n));
```
$$-(\,8\,n+9\,)\,(\,8\,n+7\,)\,\mathrm{S}(\,n+1\,) + 4\,\mathrm{S}(\,n\,)\,(\,4\,n+5\,)\,(\,4\,n+3\,) = 0\,.$$

From the resulting recurrence equation which is valid for the above sum, one can easily read off the term ratio $S(n+1)/S(n)$, hence the parameters of the hypergeometric sum. Using the package `code`, this can be automatically done by

```
> Closedform(hyperterm([1/2,-2*n],[2*n+3/2],3+2*sqrt(2),k),k,n);
```
$$\mathrm{Hyperterm}\left(\left[\frac{3}{4}, \frac{5}{4}, 1\right], \left[\frac{7}{8}, \frac{9}{8}\right], 1, n\right)\,.$$

But Koornwinder's question is a different one: He asks whether rational input (the summand of a hypergeometric sum) can generate nonrational output (a hypergeometric term representing the sum under consideration). Obviously rational input generates a rational recurrence equation. But the question remains whether the term ratio $S(n+1)/S(n)$ equivalent to this recurrence equation can be factorized rationally. All known examples in the literature are of this type. Note that rational factorization can be done algorithmically and is accessible in computer algebra systems like Maple.

## 2 Hypergeometric Product Formulas

We have discovered our counterexamples in connection with the examination of product formulas.

As an example we consider *Clausen's formula*

$$_2F_1\left(\begin{array}{c}a,b \\ a+b+1/2\end{array} \bigg| x\right)^2 = {}_3F_2\left(\begin{array}{c}2a, 2b, a+b \\ a+b+1/2, 2a+2b\end{array} \bigg| x\right). \tag{1}$$

Normally Clausen's formula is deduced by showing that both sides satisfy the same differential equation (of third order) with respect to $x$. Having Zeilberger's algorithm at hand, there is a second possibility, though. Represent the left hand side as *Cauchy product*. This is a representation by a double sum. Find a recurrence equation for the *inner sum* with respect to the summation variable of the outer sum. If the resulting recurrence equation is of first order, then it determines a hypergeometric term, the summand of the right hand sum. Hence the right hand sum constitutes a hypergeometric function whose parameters can be read off directly from the detected recurrence equation.

For the coefficient of $x^k$ of the Clausen product on the left hand side of (1), represented as Cauchy product, we get

```
> sumrecursion(
> hyperterm([a,b],[a+b+1/2],1,j)*hyperterm([a,b],[a+b+1/2],1,k-j),j,S(k));
```

$$-(k+1)(2a+k+2b)(2a+2b+1+2k)\operatorname{S}(k+1)$$
$$+2\operatorname{S}(k)(2b+k)(k+2a)(a+k+b)=0$$

This recurrence equation is valid for the coefficient of the outer sum, and hence this coefficient is given by

```
> Closedform(
> hyperterm([a,b],[a+b+1/2],1,j)*hyperterm([a,b],[a+b+1/2],1,k-j),j,k);
```

$$\operatorname{Hyperterm}\left([2b, 2a, a+b], \left[a+b+\frac{1}{2}, 2a+2b\right], 1, k\right).$$

This hypergeometric term is the coefficient of $x^k$ of the right hand sum of (1), which finishes the proof of (1).

All the usual product theorems of hypergeometric series can be deduced by this method, see [2], Chapter 7. Note that the given procedure has the advantage over the differential equation approach that it *generates* the right hand sides, given the left hand products. We will consider more examples soon where we *detect* such representations without prior knowledge.

The calculations

```
> Closedform(hyperterm([],[],x,j)*hyperterm([],[],y,k-j),j,k);
```
$$\operatorname{Hyperterm}([\,],[\,], x+y, k)$$

and

```
> Closedform(hyperterm([a],[b],x,j)*hyperterm([],[],-x,k-j),j,k);
```
$$\operatorname{Hyperterm}([-a+b],[b], -x, k)$$

generate the product identities

$$_0F_0\left(\begin{array}{c}-\\-\end{array}\bigg|\,x\right) \cdot {_0F_0}\left(\begin{array}{c}-\\-\end{array}\bigg|\,y\right) = {_0F_0}\left(\begin{array}{c}-\\-\end{array}\bigg|\,x+y\right) \qquad (2)$$

and

$$_1F_1\left(\begin{array}{c}a\\b\end{array}\bigg|\,x\right) \cdot {_0F_0}\left(\begin{array}{c}-\\-\end{array}\bigg|\,-x\right) = e^{-x} \cdot {_1F_1}\left(\begin{array}{c}a\\b\end{array}\bigg|\,x\right) = {_1F_1}\left(\begin{array}{c}b-a\\b\end{array}\bigg|\,-x\right). \qquad (3)$$

Obviously the first one is the addition theorem of the exponential function, and the second one is Kummer's identity.

The computations above deduced Equations (2) and (3) by converting

$$\sum_{j=0}^{k} \frac{x^j}{j!} \frac{y^{k-j}}{(k-j)!}$$

to $(x+y)^k/k!$, and by converting

$$\sum_{j=0}^{k} \frac{(a)_j}{(b)_j} \frac{x^j}{j!} \frac{(-x)^{k-j}}{(k-j)!}$$

to
$$\frac{(b-a)_k}{(b)_k}\frac{(-x)^k}{k!},$$

respectively.

The calculation

```
> sumrecursion(hyperterm([a],[b],x,j)*hyperterm([a],[b],-x,k-j),j,k);
```

$$(k+2)(k+1+b)(k+b)(k+2b)S(k+2) + S(k)(2a+k)x^2(2a-k-2b) = 0$$

shows that

$${}_1F_1\left(\begin{array}{c}a\\b\end{array}\bigg|\,x\right)\cdot{}_1F_1\left(\begin{array}{c}a\\b\end{array}\bigg|-x\right) = {}_2F_3\left(\begin{array}{c}b-a,a\\b,\frac{b}{2},\frac{b+1}{2}\end{array}\bigg|\,\frac{x^2}{4}\right).$$

Since the product is an even function, this can be also computed by

```
> Closedform(hyperterm([a],[b],x,j)*hyperterm([a],[b],-x,2*k-j),j,k);
```

$$\text{Hyperterm}\left([-a+b,a],\left[b,\frac{1}{2}+\frac{1}{2}b,\frac{1}{2}b\right],\frac{1}{4}x^2,k\right).$$

## 3  Counterexamples

Now we investigate examples with nonrational results. A simple class is given by

```
> Closedform(hyperterm([a],[a-1],1,j)*hyperterm([a],[a-1],1,k-j),j,k);
```

$$\text{Hyperterm}\left(\left[2a-\frac{3}{2}-\frac{1}{2}\sqrt{9-8a},2a-\frac{3}{2}+\frac{1}{2}\sqrt{9-8a}\right],\right.$$
$$\left.\left[2a-\frac{5}{2}-\frac{1}{2}\sqrt{9-8a},2a-\frac{5}{2}+\frac{1}{2}\sqrt{9-8a}\right],2,k\right),$$

in particular for $a = 3/4$

```
> Closedform(
> subs(a=3/4,hyperterm([a],[a-1],1,j)*hyperterm([a],[a-1],1,k-j)),j,k);
```

$$\text{Hyperterm}\left(\left[-\frac{1}{2}\sqrt{3},\frac{1}{2}\sqrt{3}\right],\left[-1-\frac{1}{2}\sqrt{3},-1+\frac{1}{2}\sqrt{3}\right],2,k\right)$$

generating the identities

$${}_1F_1\left(\begin{array}{c}a\\a-1\end{array}\bigg|\,x\right)^2 = {}_2F_2\left(\begin{array}{c}2a-3/2-\frac{\sqrt{9-8a}}{2},2a-3/2+\frac{\sqrt{9-8a}}{2}\\2a-5/2-\frac{\sqrt{9-8a}}{2},2a-5/2+\frac{\sqrt{9-8a}}{2}\end{array}\bigg|\,2x\right)$$

and

$${}_1F_1\left(\begin{array}{c}3/4\\-1/4\end{array}\bigg|\,x\right)^2 = {}_2F_2\left(\begin{array}{c}\frac{\sqrt{3}}{2},-\frac{\sqrt{3}}{2}\\\frac{\sqrt{3}}{2}-1,-\frac{\sqrt{3}}{2}-1\end{array}\bigg|\,2x\right)$$

respectively. The latter gives a simple counterexample for Koornwinder's question for which radicals cannot be avoided. In terms of summations, it reads as

$$\sum_{j=0}^{k}\frac{(3/4)_j}{(-1/4)_j)\,j!}\frac{(3/4)_{k-j}}{(-1/4)_{k-j}\,(k-j)!} = \frac{1-4k}{k!}\sum_{j=0}^{k}\frac{(3/4)_j\,(5/4-k)_j\,(-k)_j}{(1/4-k)_j\,(-1/4)_j}\frac{(-1)^j}{j!} =$$

$$\frac{1-4k}{k!}\,{}_3F_2\left(\begin{array}{c}3/4,5/4-k,-k\\1/4-k,-1/4\end{array}\bigg|-1\right) = \frac{(\sqrt{3}/2)_k\,(-\sqrt{3}/2)_k}{(\sqrt{3}/2-1)_k\,(-\sqrt{3}/2-1)_k}\frac{2^k}{k!},$$

and the generic example involving the parameter $a$ reads as

$$\sum_{j=0}^{k} \frac{(a)_j}{(a-1)_j\, j!} \frac{(a)_{k-j}}{(a-1)_{k-j}\,(k-j)!} = \frac{a+k-1}{(a-1)\,k!} \sum_{j=0}^{k} \frac{(a)_j\,(2-a-k)_j\,(-k)_j}{(a-1)_j\,(1-a-k)_j} \frac{(-1)^j}{j!} = \quad (4)$$

$$\frac{a+k-1}{(a-1)\,k!}\, {}_3F_2\!\left(\begin{array}{c} a, 2-a-k, -k \\ a-1, 1-a-k \end{array} \bigg| -1 \right) = \frac{\left(2a-3/2-\frac{\sqrt{9-8a}}{2}\right)_k \left(2a-3/2+\frac{\sqrt{9-8a}}{2}\right)_k}{\left(-5/2+2a-\frac{\sqrt{9-8a}}{2}\right)_k \left(-5/2+2a+\frac{\sqrt{9-8a}}{2}\right)_k} \frac{2^k}{k!}.$$

These are ${}_3F_2$-evaluations by nonrational hypergeometric terms.

Note that the hypergeometric representation given in (4) can be deduced by the Maple command

```
> simpcomb(Sumtohyper(hyperterm([a],[a-1],1,j)*hyperterm([a],[a-1],1,k-j),j));
```

$$\frac{(a+k-1)\,\mathrm{Hypergeom}([\,a,-a+2-k,-k\,],[\,a-1,-a-k+1\,],-1)}{(a-1)\,\Gamma(k+1)}$$

Similarly for $n = 2, 3, \ldots$ the products

$${}_1F_1\!\left(\begin{array}{c} a \\ a-n \end{array} \bigg| x\right)^2$$

lead to a recurrence equation of first order with respect to $k$ with polynomial coefficients of degree $2n$ in $k$ that have no proper factorization over $\mathbb{Q}$, e.g. for $n = 2$

```
> sumrecursion(hyperterm([a],[a-2],1,j)*hyperterm([a],[a-2],1,k-j),j,S(k));
```

$$-(k+1)(k^4 - 14k^3 + 8ak^3 + 67k^2 - 80ak^2 + 24a^2k^2 - 118k + 232ak + 32a^3k$$
$$- 152a^2k + 208a^2 - 96a^3 + 16a^4 - 192a + 64)S(k+1) + 2S(k)(-22k$$
$$+ 96ak - 56ak^2 - 104a^2k + 80a^2 + 31k^2 - 32a - 10k^3 - 64a^3 + 16a^4$$
$$+ 32a^3k + 24a^2k^2 + k^4 + 8ak^3) = 0$$

or for $n = 3$

```
> sumrecursion(hyperterm([a],[a-3],1,j)*hyperterm([a],[a-3],1,k-j),j,S(k));
```

$$-(k+1)(-264k^4 a - 8448a - 5136k + 289k^4 + 6180a^2k^2 + 2180ak^3 - 9216a^3$$
$$- 2016a^3k^2 - 1032a^2k^3 + 7808a^3k + 3712a^4 - 27k^5 + 2304 + 12k^5 a + k^6$$
$$- 1968a^4 k - 768a^5 + 13984ak + 4102k^2 + 192a^5 k + 240a^4 k^2 + 160a^3 k^3$$
$$+ 60k^4 a^2 - 15000a^2 k + 12352a^2 - 8232ak^2 - 1533k^3 + 64a^6)S(k+1) + 2$$
$$\mathrm{S}(k)(-204k^4 a - 768a - 504k + 169k^4 + 3444a^2 k^2 + 1244ak^3 - 3264a^3$$
$$- 1536a^3 k^2 - 792a^2 k^3 + 4256a^3 k + 1984a^4 - 21k^5 + 12k^5 a + k^6 - 1488a^4 k$$
$$- 576a^5 + 3064ak + 982k^2 + 192a^5 k + 240a^4 k^2 + 160a^3 k^3 + 60k^4 a^2$$
$$- 5496a^2 k + 2560a^2 - 3156ak^2 - 627k^3 + 64a^6) = 0$$

To present some other results of similar type, we calculate the product

$${}_2F_1\!\left(\begin{array}{c} a, a+b \\ a-2 \end{array} \bigg| x\right) \cdot {}_2F_1\!\left(\begin{array}{c} a, a+b \\ a-2 \end{array} \bigg| -x\right)$$

which is representable by a hypergeometric function with parameters using square roots, but in a rather complicated way:

```
> sumrecursion(
> hyperterm([a,a+b],[a-2],x,j)*hyperterm([a,a+b],[a-2],-x,2*k-j),j,S(k));
```

$$-(k+1)(-7ak^2 + 5k^2b + k^2b^2 + a^2k^2 + 8k^2 - 2abk^2 + 21ak + 4abk - 12k$$
$$- 2ka^2b - 13a^2k + 2ka^3 - 5kb - kb^2 + 4 - 12a + a^4 - 6a^3 + 13a^2)S(k+1)$$
$$+ S(k)(b+k+a)(a^4 - 4a^3 + 2ka^3 - 11a^2k + a^2k^2 - 2a^2b + a^2 - 2ka^2b$$
$$+ 2a - 7ak^2 + 7ak + 2ab - 2abk^2 + 8k^2 + 5kb + k^2b^2 + kb^2 + 5k^2b + 4k)$$
$$x^2 = 0$$

generating the hypergeometric term

```
> factor(Closedform(
> hyperterm([a,a+b],[a-2],x,j)*hyperterm([a,a+b],[a-2],-x,2*k-j),j,k));
```

$$\text{Hyperterm}\left(\left[a+b, \frac{1}{2}\frac{2a^3 - 11a^2 - 2a^2b + 7a + 5b + b^2 + 4 - B}{A},\right.\right.$$
$$\left.\frac{1}{2}\frac{2a^3 - 11a^2 - 2a^2b + 7a + 5b + b^2 + 4 + B}{A}\right], \left[\right.$$
$$\frac{1}{2}\frac{-13a^2 + 4ab + 21a - 12 - 2a^2b + 2a^3 - b^2 - 5b - B}{A},$$
$$\left.\frac{1}{2}\frac{-13a^2 + 4ab + 21a - 12 - 2a^2b + 2a^3 - b^2 - 5b + B}{A}\right], x^2, k\right)$$

with
$$A := b^2 + 8 + 5b - 2ab + a^2 - 7a$$

and

$$B := \left((a+1+b)(a^3 + 9a^2 + 4a^2b^2 + 15a^2b - 24a - 34ab - 9ab^2 + 24b + 9b^2 + b^3 + 16)\right)^{1/2}.$$

Similarly

```
> sumrecursion(
> hyperterm([a,a+b],[a-3],1,j)*hyperterm([a,a+b],[a-3],-1,2*k-j),j,S(k));
```

$$-(k+1)(36 - 132a - 132k + 58a^4 - 84abk^2 - 144a^3 + a^6 - 54ka^2b + 66abk$$
$$+ 168a^2k^2 + 3ab^2k^3 + 3a^2k^2b^2 + 193a^2 - 9k^3b^2 + 42a^2bk^2 + 24abk^3$$
$$+ 18a^3bk - 3a^2bk^3 - 6a^3bk^2 - 15a^2k^3 - 42a^3k^2 - 39a^4k + a^3k^3 + 3a^4k^2$$
$$+ 56ak^3 + 144k^2 - 351a^2k + 81k^2b + 3a^5k - 3a^4bk - 12a^5 - 2kb^3$$
$$- 255ak^2 + 170ka^3 + 3k^2b^3 - 9ak^2b^2 - 3a^2kb^2 - k^3b^3 - 48k^3 + 6ab^2k$$
$$- 15kb^2 + 24k^2b^2 - 49kb + 337ak - 32k^3b)S(k+1) + S(k)(b+k+a)(6a$$
$$+ 12k + 22a^4 - 12abk^2 - 15a^3 + a^6 + 21ka^2b - 30abk + 123a^2k^2$$
$$+ 3ab^2k^3 + 3a^2k^2b^2 - 5a^2 - 9k^3b^2 + 33a^2bk^2 + 24abk^3 + 6a^3bk$$
$$- 3a^2bk^3 - 6a^3bk^2 - 15a^2k^3 - 39a^3k^2 - 33a^4k + a^3k^3 + 3a^4k^2 + 56ak^3$$
$$- 15a^2b - 60a^2k - 15k^2b + 3a^5k - 3a^4bk - 9a^5 + kb^3 + 6ab - 87ak^2$$
$$+ 89ka^3 + 3a^2kb^2 - k^3b^3 - 48k^3 - 3a^4b - 3ab^2k + 6kb^2 - 3k^2b^2 + 17kb$$
$$- 5ak - 32k^3b + 12a^3b) = 0,$$

and first order recurrence equations, more and more complicated with increasing $n$, can be determined iteratively for the coefficients of

$$_2F_1\left(\begin{array}{c}a, a+b\\a-n\end{array}\bigg|\, x\right) \cdot {}_2F_1\left(\begin{array}{c}a, a+b\\a-n\end{array}\bigg|\, -x\right).$$

We finish this note with a rather simple family of examples. If $n \in \mathbb{N}$ is a nonnegative integer, then the polynomials ($k \in \mathbb{N}$)

$$_2F_1\left(\begin{array}{c}a+n, -k\\a\end{array}\bigg|\, x\right)$$

satisfy a recurrence equation of first order with respect to $k$ with polynomial coefficients of degree $n$ in $k$ that have no proper factorization over $\mathbb{Q}$. As a particular case one has for $n = 3$

```
> sumrecursion(hyperterm([a+3, -k],[a],x,j),j,S(k));
```

$(-2\,a - 6\,x^2\,k\,a^2 + 6\,x^3\,k\,a + 3\,x^3\,k\,a^2 + 2\,x^3\,k + 2\,x^3\,a - 3\,x^2\,a^3 + x^3\,k^3 + 3\,x^3\,k^2$
$\quad - 15\,x^2\,k\,a - 3\,x^2\,k^2\,a + 3\,x^3\,k^2\,a - 6\,x^2\,k^2 + 3\,x^3\,a^2 + x^3\,a^3 + 9\,x\,k\,a + 6\,x\,a + 6\,x\,k$
$\quad - 3\,a^2 - a^3 + 3\,x\,k\,a^2 + 3\,x\,a^3 + 9\,x\,a^2 - 6\,x^2\,a - 6\,x^2\,k - 9\,x^2\,a^2)\mathrm{S}(\,k+1\,) + \mathrm{S}(\,k\,)$
$(\,x-1\,)(-2\,a + 6\,x + 6\,x^3 - 6\,x^2\,k\,a^2 + 12\,x^3\,k\,a + 3\,x^3\,k\,a^2 + 11\,x^3\,k + 11\,x^3\,a$
$\quad - 3\,x^2\,a^3 + x^3\,k^3 + 6\,x^3\,k^2 - 21\,x^2\,k\,a - 3\,x^2\,k^2\,a + 3\,x^3\,k^2\,a - 6\,x^2\,k^2 + 6\,x^3\,a^2$
$\quad + x^3\,a^3 + 9\,x\,k\,a + 15\,x\,a + 6\,x\,k - 3\,a^2 - a^3 + 3\,x\,k\,a^2 - 12\,x^2 + 3\,x\,a^3 + 12\,x\,a^2$
$\quad - 24\,x^2\,a - 18\,x^2\,k - 15\,x^2\,a^2) = 0$